\documentclass{amsart} 

\usepackage{amssymb} 

\newtheorem{theorem}{Theorem}[section] 
\newtheorem{lemma}[theorem]{Lemma} 
\newtheorem{proposition}[theorem]{Proposition} 
\newtheorem{corollary}[theorem]{Corollary} 
\theoremstyle{definition} 
\newtheorem{definition}[theorem]{Definition} 
\newtheorem{example}[theorem]{Example} 
 
\theoremstyle{remark} 
\newtheorem{remark}[theorem]{Remark} 
\numberwithin{equation}{section} 

\def\cal#1{\mathcal{#1}}

\def\frk{\frak}               

\def\pp{{\frk p}}

\def\mm{{\frk m}}

\def\nn{{\frk n}}

\def\opn#1#2{\def#1{\operatorname{#2}}} 
\opn\chara{char}
\opn\length{\ell}
\opn\pd{pd}
\opn\rk{rk}
\opn\projdim{proj\,dim}
\opn\rank{rank}
\opn\depth{depth}
\opn\grade{grade}
\opn\height{ht}
\opn\embdim{emb\,dim}
\opn\codim{codim}

\opn\Tr{Tr}
\opn\bigrank{big\,rank}
\opn\superheight{superheight}\opn\lcm{lcm}
\opn\trdeg{tr\,deg}%
\opn\reg{reg}
\opn\lreg{lreg}
\opn\div{div}
\opn\Div{Div}
\opn\WDiv{WDiv}
\opn\cl{cl}
\opn\Cl{Cl}
\opn\Spec{Spec}
\opn\Supp{Supp}
\opn\supp{supp}
\opn\Sing{Sing}
\opn\Ass{Ass}
\opn\Ann{Ann}
\opn\Rad{Rad}
\opn\Soc{Soc}
\opn\Ker{Ker}
\opn\Coker{Coker}
\opn\Im{Im}
\opn\Hom{Hom}
\opn\Tor{Tor}
\opn\Ext{Ext}
\opn\End{End}
\opn\Aut{Aut}
\opn\id{id}

\opn\nat{nat}
\opn\pff{pf}
\opn\Pf{Pf}
\opn\GL{GL}
\opn\SL{SL}
\opn\mod{mod}
\opn\ord{ord}
\opn\Proj{Proj}
\opn\aff{aff}
\opn\con{conv}
\opn\relint{relint}
\opn\st{st}
\opn\lk{lk}
\opn\cn{cn}
\opn\core{core}
\opn\vol{vol}
\opn\link{link}
\opn\star{star}
\opn\gr{gr}

\def\pot#1#2{#1[\kern-0.28ex[#2]\kern-0.28ex]}

\opn\dirlim{\underrightarrow{\lim}}
\opn\inivlim{\underleftarrow{\lim}}
\let\iso=\cong
\let\Union=\bigcup

\let\mcone= * 
\def\Implies{\ifmmode\Longrightarrow \else
     \unskip${}\Longrightarrow{}$\ignorespaces\fi}
\def\implies{\ifmmode\Rightarrow \else
     \unskip${}\Rightarrow{}$\ignorespaces\fi}
\def\iff{\ifmmode\Longleftrightarrow \else
     \unskip${}\Longleftrightarrow{}$\ignorespaces\fi}

\let\:=\colon
\opn\H{H}
\opn\Pic{Pic}

\let\epsilon\varepsilon
\opn\inii{in}
\opn\inim{inm}
\opn\set{set}
\def\pnt{{\raise0.5mm\hbox{\large\bf.}}}

\def\Coh#1#2{H_{\mm}^{#1}(#2)}

\begin{document} 

\title{Stanley-Reisner ideals whose powers have finite length cohomologies}
\author{Shiro Goto} 
\address{Department of Mathematics, 
         School of Science and Technology,
         Meiji University, 
         214-8571, Japan} 
\email{goto@math.meiji.ac.jp} 
\author{Yukihide Takayama} 
\address{Department of Mathematical Sciences, 
         Ritsumeikan University, 
         1-1-1 Nojihigashi, Kusatsu, Shiga 525-8577, Japan} 
\email{takayama@se.ritsumei.ac.jp}



\begin{abstract} 
We introduce a class of Stanley-Reisner ideals called generalized 
complete intersection, which is characterized by the property that 
all the residue class rings of powers of the ideal have FLC. 
We also give a combinatorial characterization of such ideals.  
\end{abstract}

\maketitle 


\section{Introduction}

Let $S = K[X_1,\ldots, X_n]$ be a polynomial ring over a field $K$ and
$\mm = (X_1,\ldots, X_n)$. We view $(S, \mm)$ as a standard graded
algebra with the unique graded maximal ideal $\mm$.  Let $I\subset S$
be a graded ideal and set $R = S/I$.  Then $R$, or $I$, is called generalized
Cohen-Macaulay or FLC if the local cohomology $\Coh{i}{R}$ has 
finite length for all $i<d=\dim R$. If $I$ is a Stanley-Reisner ideal,
i.e., a square-free monomial ideal, it is well known that FLC 
coincides with Buchsbaumness, and Buchsbaum Stanley-Reisner ideals are well
understood through topological characterization of the simplicial
complexes corresponding  the ideal. However, FLC monomial ideals,
which are not always square-free, have not been well understood. In
\cite{T}, the second author gave combinatorial characterizations of
FLC monomial ideals for $d\leq 3$ and a method for constructing FLC
monomial ideals from Buchsbaum Stanley-Reisner ideals. But the problem
to find  fairly large classes of FLC monimial ideals has been open.

The aim of this paper is to give an answer to this problem.  As shown
in \cite{T}, for a monomial ideal $I\subset S$, it is FLC only when
$\sqrt{I}$ is a Buchsbaum Stanley-Reisner ideal. One of the typical
cases of this situation is that $J\subset S$ is a Buchsbaum
Stanley-Reisner ideal and we consider its powers, i.e., $I = J^\ell$ for
$\ell = 1,2,\ldots$. Unfortunately, such monomial ideal $J$ is not always FLC.
We show that $I$ is FLC for all integer $\ell\geq 1$ if and only if the
simplicial complex $\Delta$ corresponding to $J$ is pure and
$k[\lk_{\Delta}(\{i\})]$ is complete intersection for all
$i\in[n]=\{1,\ldots, n\}$ (Theorem~\ref{theorem:main1}). Then 
we classify such simplicial complexes $\Delta$ (Theorem~\ref{theorem:main2}).

For a Noetherian local ring $(R,\nn)$ and a finitely generated 
$R$-module $M$,  we denote by $\mathrm{Min}_R{M}$ the set of minimal
primes in $\Ass_R(M)$.  If
$M$ is an Artinial $R$-module, we denote its length as an $R$-module
by $\ell_R(M)$.

We thank J.\ Herzog who read early version of this manuscript and gave
us valuable comments.

\section{Powers of generically complete intersection ideal}

Let $I$ be an ideal of a Cohen-Macaulay local ring $(S, \nn)$ such
that $S/I$ is generically a complete intersection, namely for
arbitrary $\pp\in\mathrm{Min}_S{S/I}$, the ideal $IS_{\pp}$ is
generated by $S_{\pp}$-regular sequence. We consider here the
condition that when $S/I^\ell$, $\ell=1,2,\ldots$, are Cohen-Macaulay.
Now we recall the following well-known result,
which was originally proved  by Cowsik and Nori
\cite{CN} and whose refinement as stated in the following lemma 
has been given, for example, in \cite{AV, HO, Wal}.

\begin{lemma}
Let $I\subset S$ be as above. Then
the following are equivalent.
\label{mainlemma}
\begin{enumerate}
\item[$(1)$] $I$ is generated by an $S$-regular sequence,
\item[$(2)$] $S/I^{\ell + 1}$ is Cohen-Macaulay
for arbitrary integer $\ell \geq 0$,
\item[$(3)$] $\Lambda 
= \{\ell \geq 0 
   \mid S/I^{\ell + 1} \text{ is Cohen-Macaulay} \}$
is an infinite set.
\end{enumerate}
\end{lemma}

Notice that Stanley-Reisner ideal is a typical class of generically
complete intersection (cf. Th.~5.1.4 \cite{BH}).  
Let $\Delta$ be a simplicial complex over the vertex set 
$[n]=\{1,\ldots, n\}$. We always assume 
that $\{i\}\in\Delta$ for all $i\in [n]$
unless otherwise stated.
Then we denote the Stanley-Reisner ring 
corresponding to $\Delta$ by $k[\Delta] = S/I_\Delta$,
where $S = k[X_1,\ldots, X_n]$ is the polynomial ring over the field $k$
and $I_\Delta\subset S$ is the square-free monomial ideal 
corresponding to $\Delta$.
For a monomial ideal $I\subset S$, we denote by $G(I)$ the minimal 
set of generators.
Also for $F\in\Delta$, 
we define the {\em link} by $\lk_{\Delta}(F) = \{G \;\vert\; G\cup F\in\Delta, G\cap F=\emptyset\}$.
See \cite{BH} Chapter~5 for the detail on these terminologies.

Now we introduce the following notion.
\def\gci{generalized complete intersection }
\begin{definition}
A Stanley-Reisner ring $k[\Delta]$ is 
called a {\em \gci} $($gCI$)$ if $\Delta$ is pure 
and $k[\lk_{\Delta}(\{i\})]$ is complete intersection 
for all $i\in[n]$.
We also call $\Delta$ to be generalized Cohen-Macaulay
if $k[\Delta]$ is gCI.
\end{definition}

This terminology comes from the analogy of generalized
Cohen-Macaulayness.  We notice the following fact, whose proof is left
to the readers.

\begin{proposition}
\label{gci}
For a simplicial complex $\Delta$, the followings
are equivalent:
\begin{enumerate}
\item [$(i)$] $\Delta$ is  gCI.
\item [$(ii)$] $\Delta$ is pure and $k[\Delta]_P$ is 
a complete intersection for every prime $P(\ne\mm)$.
\item [$(iii)$] $\Delta$ is pure and $k[\Delta]_{X_i}$ is 
a complete intersection for every $i=1,\ldots, n$.
\item [$(iv)$] $\Delta$ is pure and 
if $G(I_\Delta)  = \{u_1,\ldots, u_\ell\}$ then
$\bar{u}_1,\ldots, \bar{u}_\ell$ is a $S_i$-regular
sequence for all $i\in[n]$,
where $S_i = k[X_1,\ldots,\overset{i}{\vee}\ldots, X_n]$
and $\bar{u}_j$ is obtained from $u_j$ by substituting $X_i$ by $1$.
\end{enumerate}
\end{proposition}

\begin{corollary}
\label{ci-is-genci}
If $k[\Delta]$ is  complete intersection, then it is
 \gci.
\end{corollary}

Now we can prove the main theorem of this section.

\begin{theorem}
\label{theorem:main1}
Let $n \geq 1$ be an integer
and let $\Delta$ be a simplicial complex
over the vertex set $[n]$. Then following are equivalent.
\begin{enumerate}
\item[$(1)$] $k[\Delta]$ is \gci
\item[$(2)$] 
$S/{I_{\Delta}}^{\ell + 1}$ has FLC for arbitrary integer $\ell \geq 0$
\item[$(3)$] The set 
$\{ \ell \geq 0 \mid \text{ the ring $S/{I_{\Delta}}^{\ell + 1}$
 has FLC} \}$ is infinite.
\end{enumerate}
If one of these conditions holds, 
$k[\Delta]$ is Buchsbaum.
\end{theorem}

\begin{proof}
\begin{description}
\item [$(1)\Rightarrow (2)$] 
Since $k[\Delta]_{x_i} \iso k[x_i, x_i^{-1}][\lk(\{i\})]$, 
$(1)$ is equivalent to the condition that 
$\Delta$ is pure and $k[\Delta]_P$ is complete intersection
for all prime ideals $P\ne\mm$. In particular, this implies that 
$k[\Delta]$ is Buchsbaum. See, for example, Excer.~5.3.6~\cite{BH}.
Also by Lemma~\ref{mainlemma}, $S_P/(I_{\Delta}S_P)^\ell$ is 
Cohen-Macaulay for all integers  $\ell\geq 0$. This implies that 
$S/I_\Delta^\ell$ is generalized Cohen-Macaulay, i.e., it has FLC.
\item [$(3)\Rightarrow (1)$]
$(3)$ implies that, for arbitrary $P\ne\mm$,
$S_P/(I_{\Delta}S_P)^\ell$ is Cohen-Macaulay for infinitely
may integers $\ell>0$. Thus by Lemma~\ref{mainlemma}
$k[\Delta]_P = S_P/I_{\Delta}S_P$ must be complete intersection, and 
this implies that $k[\Delta]_{x_i} \iso k[x_i,x_i^{-1}][\lk(\{x_i\})]$
is complete intersection for all $i=1,\ldots, n$. Then we 
obtain $(1)$.
\end{description}
\end{proof}

\section{Combinatorial characterization}

The aim of this section is to give a combinatorial
characterization of gCI simplicial complexes.
We denote by $\Delta$ a simplicial complex over the vertex
set $[n]$ and $S = k[X_1,\ldots, X_n]$ is the polynomial ring
over the field $k$.

\subsection{core of simplicial complex}
We recall here the notion of {\em core} of simplicial complexes
(\cite{BH}~section~5.5).
For $F\in\Delta$, we define $\star_{\Delta}(F) = \{G \;\vert\;
G\cup F\in\Delta \}$.
We also define $\core [n] = \{i\in [n]\;\vert\; \star_\Delta(i)\ne \Delta\}$.
Then the core of $\Delta$ is defined by 
$\core\Delta = \{F\cap \core [n] \;\vert\; F\in\Delta\}$.

For two simplicial complexes $\Delta_1$ and $\Delta_2$,
we define the {\em join} of them by 
$\Delta_1 * \Delta_2
= \{ F\cup G \;\vert\;  F\in\Delta_1,\; G\in \Delta_2\}$.
In particular, if $\Delta_1 = \{\emptyset, \{i\}\}$, which 
we will abvreviate simply to $\Delta_1 = \{i\}$, $\{i\}*\Delta$
is called a {\em cone} over $\Delta$.
Notice that if $i\in[n]\backslash\core[n]$ then we have 
\begin{equation*}
\Delta = \star_\Delta(i) = \{G\;\vert\; G\cup\{i\}\in\Delta\}
       = \{i\}*\link_\Delta(i).
\end{equation*}
Thus if $\Delta\ne \core\Delta$
then $U:=[n]\backslash\core[n]$ is non-empty
and any element from $G(I)$ does not contain $X_j$, $j\in U$.
This implies that $k[\Delta]$ is a polynomial ring over the ring 
$k[\core\Delta]$:
$k[\Delta] = k[\core\Delta][X_j \;\vert\; j\in U]$
with $k[\core\Delta] = S'/I_\Delta S'$, $S' = k[X_i\;\vert\; i\in\core[n]]$.

\subsection{complete intersection by localization}

Let $I_\Delta = (u_1,\ldots, u_\ell) \subset S$ be a
square-free \gci ideal that is not a complete intersection.
Then, $(a)$ $\Delta$ is pure, and  $(b)$
the localization by $X_i$, $i\in [n]$, is complete intersection.
In this subsection, we consider the problem to characterize 
$G(I_\Delta)$ with the ideal $I_\Delta$ having the property $(b)$.
The purity condition $(a)$ will be considered in the next subsection.
First of all we have
\begin{lemma}
\label{lemma0}
We must have $\ell \geq 2$ and $\Delta = \core\Delta$.
\end{lemma}
\begin{proof} The necessity of the first condition is clear.
Now assume that 
$\Delta\ne\core\Delta$. Then there exists an element $i\in
[n]\backslash\Union_{j=1}^\ell\supp(u_j)$.  $k[\Delta]_{X_i}$ must
be a complete intersection and $u_1,\ldots, u_\ell$ remain to be 
the minimal set of generators in $I_\Delta S_{X_i}$
so that $\supp(u_j)$, $1\leq j\leq \ell$, must be pairwise disjoint.
Thus 
we know that $k[\Delta]$ itself is a complete
intersection, which contradicts the assumption.
\end{proof}

Now we will assume $\Delta = \core\Delta$ in the following.  We will
rephrase our problem by purely combinatorial setting. 
Let $S_j:= \supp(u_j)$, $j=1,\ldots, \ell$ and set
\begin{equation*}
{\cal F}_{\Delta}:= \{S_1,\ldots, S_\ell\}.
\end{equation*}
Since $\Delta = \core\Delta$, we have 
\begin{equation}
\label{cond-1}
   S_1\cup\cdots\cup S_\ell = [n].
\end{equation}
Also, 
since $\{i\}\in\Delta$ for all $i\in[n]$, we have 
\begin{equation}
\label{cond0}
\sharp S_j \geq 2\quad\mbox{for all $1\leq j\leq \ell$}.
\end{equation}
Since $\{u_1,\ldots, u_\ell\} = G(I_\Delta)$, we have 
\begin{equation}
\label{cond1}
S_i\not\subset S_j\quad\mbox{for all } i,j \mbox{ with } i\ne j.
\end{equation}
Now for $i\in[n]$, we set ${\cal F}_i = \{S\backslash\{i\}\;\vert\;
S\in{\cal F}_{\Delta}\}$.  We will say that $S(\in{\cal F}_\Delta)$ is minimal in
${\cal F}_i$ if $(\emptyset\ne)S\backslash\{i\}$ is minimal in ${\cal
F}_i$ with regard to the set inclusion. We set $\min{\cal F}_i:=
\{S_j\backslash\{i\}\;\vert\; S_j\mbox{ is minimal in } {\cal F}_i\}$,
which represents $G(I_\Delta S_{X_i})$.
Since $S_{X_i}/(u_1,\ldots, u_\ell)S_{X_i}$ is a 
complete intersection for all
$i\in[n]$, we must have 
\begin{equation}
\label{cond2}
     (S\backslash\{i\})\cap (S'\backslash\{i\}) =\emptyset
\quad \mbox{for all distinct }
S\backslash\{i\}, S'\backslash\{i\}\in\min{\cal F}_i
\quad\mbox{for all $i\in[n]$}.
\end{equation}
Finally, since we do not assume that $k[\Delta]$ itself is 
a complete intersection, we exclude the case that 
elements in ${\cal F}_{\Delta}$ are pairwise disjoint so that  we pose
\begin{equation}
\label{cond3}
 S\cap S' \ne\emptyset\quad\mbox{for some }S, S'\in{\cal F}_{\Delta}.
\end{equation}

Now our problem is as follows: determine ${\cal F}_{\Delta} = \{S_1,\ldots,
S_\ell\}$, $S_j\subset [n]$, satisfying the conditons 
$(\ref{cond-1})$, $(\ref{cond0})$,
$(\ref{cond1})$, $(\ref{cond2})$ and $(\ref{cond3})$, together with
some condition assuring the purity of  $\Delta$.

\begin{lemma}
\label{lemma1}
For every element $S\in{\cal F}_{\Delta}$, there exists another element
$S'\in{\cal F}_{\Delta}$ such that $S\cap S'\ne\emptyset$.
\end{lemma}
\begin{proof}
Assume that $S_1\cap S_j=\emptyset$ 
for $j=2,\ldots, \ell$.
Then for any $i\in S_1$ we have $i\notin S_j$
for $j=2,\ldots, \ell$.
Thus ${\cal F}_i = \{S_1\backslash \{i\}, S_2,\ldots, S_\ell \}$
and, by the condition $(\ref{cond1})$ and the assumption on $S_1$ we know 
that ${\cal F}_i = \min {\cal F}_i$.
Thus, by the condition $(\ref{cond2})$,
we know that $S_1,\ldots, S_\ell$ are pairwise disjoint, which 
contradicts the condition $(\ref{cond3})$. 
\end{proof}

\begin{lemma}
\label{lemma2}
Let $S\in{\cal F}_{\Delta}$ be arbitrary. For any $i\in S$,
$S$ is minimal in ${\cal F}_i$.
Also if $S\backslash\{i\}\subset S'\backslash\{i\}$
for distinct $S, S'\in{\cal F}_\Delta$, we have $i\in S$.
\end{lemma}
\begin{proof}
Assume that $S$ with $i\in S$ is not minimal in ${\cal F}_i$.
Then there exists another element
$S'\in{\cal F}_{\Delta}$ such that $S'\backslash\{i\}\subset S\backslash\{i\}$.
We must have $i\in S'$ and $i\notin S$, since otherwise we have
$S'\subset S$ which contradicts $(\ref{cond1})$. But this is impossible
since we have $i\in S$ by assumption.
Thus $S$ must be minimal
in ${\cal F}_i$ with $i\in S$. 
\end{proof}

\begin{lemma}
\label{lemma3}
For every $S, S'\in{\cal F}_{\Delta}$, $\sharp(S\cap S')\leq 1$.
\end{lemma}
\begin{proof}
Assume that $i, j\in S_1\cap S_2$, $i\ne j$.
Since $j\in (S_1\backslash\{i\})\cap (S_2\backslash\{i\})\ne\emptyset$,
either $S_1$ or $S_2$ is non-minimal in ${\cal F}_i$ by the condition
$(\ref{cond2})$.  But, since $i\in S_1$ and $i\in S_2$,
they must be minimal in ${\cal F}_i$ by Lemma~\ref{lemma2}, a contradiction.
Thus $S_1\cap S_2$ contains at most one element.
\end{proof}

\begin{lemma}
\label{lemma4}
Let $S\in{\cal F}_{\Delta}$ be such that $\sharp{S}\geq 3$.
If $S'\in{\cal F}_{\Delta}$ is another element such that 
$S'\cap S\ne\emptyset$, then $\sharp{S'}=2$.
\end{lemma}
\begin{proof}
Let $S_1, S_2\in {\cal F}_{\Delta}$ be such that 
$\sharp{S_1}\geq 3$ and $S_1\cap S_2\ne\emptyset$.
Notice that the existence of such $S_2$ is assured
by Lemma~\ref{lemma1}. 
Now by Lemma~\ref{lemma3}, we can assume
without loss of generality that $S_1\cap S_2=\{1\}$ and 
$2,3\in S_1\backslash S_2$. 
%
First we consider ${\cal F}_2$.
Since $(S_1\backslash\{2\})\cap
(S_2\backslash\{2\})=\{1\}\ne\emptyset$
and $2\in S_1$, we must have
$(S_3\backslash\{2\})\subset (S_2\backslash\{2\})$
for some  $S_3 (\ne S_2)\in{\cal F}_{\Delta}$ with 
$2\in S_3$ by $(\ref{cond2})$ and Lemma~\ref{lemma2}.
We know that $S_3\ne S_1$. 
By Lemma~\ref{lemma3} we know 
that $S_3 = \{2, k\}$ with some $k\in S_2\backslash S_1$. 
We can assume without loss of generality that $k=4$: $S_3 = \{2,4\}$.
%
Next we consider  ${\cal F}_1$. 
We have $(S_1\backslash\{1\})\cap (S_3\backslash\{1\})=\{2\}
\ne\emptyset$ and $1\in S_1$, so that 
by $(\ref{cond2})$ and Lemma~\ref{lemma2} 
we must have $(S_4\backslash\{1\})\subset (S_3\backslash\{1\})$
for some $S_4(\ne S_3)\in{\cal F}_{\Delta}$ with $1\in S_4$.
Thus by Lemma~\ref{lemma3} we know that the only possible 
case is 
$S_4 = S_2 = \{1,4\}$,
i.e., $\sharp S_2 =2$ as required.
\end{proof}

\begin{lemma}
\label{lemma5}
Let $S\in{\cal F}_{\Delta}$ be such that $\sharp{S}\geq 3$.
If $S' = \{i,j\}\in{\cal F}_{\Delta}$ is another element such that 
$S'\cap S = \{j\}$, then
$\{i, k\}\in {\cal F}_{\Delta}$ for all $k\in S$.
\end{lemma}
\begin{proof}
As in the proof of Lemma~\ref{lemma3}, we can assume without loss
of generality that $\sharp S_1\geq 3$, $S_2 = \{1,2\}$, 
$3,4\in S_1\backslash S_2$ and $S_1\cap S_2=\{1\}$. It suffices to 
show that $\{2,3\}\in {\cal F}_{\Delta}$.  We consider ${\cal F}_3$.
Since $3\in S_1$, $S_1$ is minimal in ${\cal F}_3$ by
Lemma~\ref{lemma2}.
We have 
$(S_1\backslash\{3\})\cap(S_2\backslash\{3\})=\{1\}\ne\emptyset$
and 
$(S_1\backslash\{3\})\not\subset(S_2\backslash\{3\})$.
Thus, by $(\ref{cond2})$ and Lemma~\ref{lemma2}, we have 
$(S_3\backslash\{3\})\subset (S_2\backslash\{3\}) = S_2$
for some $S_3 \in{\cal F}_{\Delta}$ with $3\in S_3$. 
By Lemma~\ref{lemma3},
we know that $S_3 = \{2, 3\}$.
\end{proof}

Now for an element $S\in{\cal F}_{\Delta}$ such that $\sharp{S}\geq 3$,
we set
\begin{equation*}
{\cal C}(S) = \{i\in [n]\;\vert\; \{i,j\}\in{\cal F}_{\Delta}\mbox{ for
some $j\in S$}\}.
\end{equation*}
By the condition $(\ref{cond1})$, we know ${\cal C}(S)\cap S=\emptyset$.
According to  Lemma~\ref{lemma4}, 
for any $i\in {\cal C}(S)$, we have $\{i,k\}\in{\cal F}_{\Delta}$ for all $k\in S$.

\begin{lemma}
\label{lemma6}
Let $S\in{\cal F}_{\Delta}$ be such that $\sharp{S}\geq 3$.
For an  $i\in [n]$, if 
$i\notin{\cal C}(S)\cup S$, then 
$\{i,k\}\in{\cal F}_{\Delta}$ for all $k\in {\cal C}(S)$.
\end{lemma}
\begin{proof}
Assume that for an $i\notin{\cal C}(S)$  there exists 
$k\in{\cal C}(S)$ such that $\{i,k\}\notin{\cal F}_{\Delta}$.
We will deduce a contradiction. 
By $(\ref{cond-1})$ there exists at least one $T\in{\cal F}_{\Delta}$ such 
that $i\in T$. Now we have two cases:
$(i)$ $\{i,k\}\subset T$ and $\sharp T\geq 3$,
or $(ii)$ for all such $T$ we have $\{i,k\}\not\subset T$.

We first consider the case $(i)$. 
We have $\{k,j\}\in{\cal F}_{\Delta}$ for some  $j\in S$,
and $k \in T$ with $\sharp{T}\geq 3$. Thus by Lemma~\ref{lemma4}
$\{i,j\}\in {\cal F}_{\Delta}$, which contradicts the assumption that 
$i\notin {\cal C}(S)$.

Next we consider the case $(ii)$. 
For any distinct $j_1, j_2\in S$, we have 
$(\{k,j_1\}\backslash\{i\})\cap(\{k,j_2\}\backslash\{i\})=\{k\}\ne\emptyset$
in ${\cal F}_i$. Thus there exists $S'\in{\cal F}_{\Delta}$ with 
$i\in S'$ such that $(S'\backslash\{i\})\subset(\{k,j_p\}\backslash\{i\})$
with $p=1$ or $2$. Thus we must have 
either $S' = \{i, k\}$ or $S' = \{i, j_p\}$. But both of the 
cases contradicts the assumptions of non-existence of $\{i,k\}\in{\cal F}_{\Delta}$
and $i\notin{\cal C}(S)$.
\end{proof}

In the following, we partly use the language of graph theory.
We call an element $E\in{\cal F}_\Delta$ with $\sharp{E}=2$ 
an {\em edge}. Also we will call a set $P = \{E_1,\ldots, E_q\}$
of edges with $E_i\cap E_{i+1}\ne\emptyset$, $i=1,\ldots, q-1$, 
a {\em path}.

\begin{lemma}
\label{lemma7}
Any two elements  $i, j\in [n]$ 
are 
linked with a path $P = \{E_1,\ldots, E_q\}$, 
such that $i\in E_1$ and $j\in E_q$.
\end{lemma}
\begin{proof} The case that $\{i,j\}\in {\cal F}_\Delta$ 
is trivial.
Also, if $\{i,j\}\subset S$ for some $S\in {\cal F}_\Delta$ with 
$\sharp{S}\geq 3$, then by Lemma~\ref{lemma1},~\ref{lemma4} and
~\ref{lemma5}, there exists $k\in {\cal C}(S)$ and 
we have $\{i,k\}, \{k,j\}\in{\cal F}_\Delta$. Namely, $i$ and $j$
are linked with the path $P = \{\{i,k\}, \{k,j\}\}$.
Now we will consider other cases.

Let $i, j\in [n]$ be arbitrary. Assume that there exists a set $Q =
\{S_1,\ldots, S_r\}$, $r\geq 2$, of elements from ${\cal F}_{\Delta}$
such that $i\in S_1$, $j\in S_r$ and $S_i\cap S_{i+1}\ne\emptyset$ for
$i=1,\ldots, r-1$, where we have $\sharp{S_j}\geq 3$ for some $1\leq
j\leq r$.
By Lemma~\ref{lemma4}, we must have $\sharp{S_{j-1}}=\sharp{S_{j+1}}=2$
if $S_{j-1}$ and $S_{j+1}$ exist. We consider the case that $S_{j-1}$
exists. Another case is similar. We set $S_{j-1} = \{i_1,i_2\}$
and $S_j = \{i_2,\ldots, i_q\}$ with $q\geq 4$. Then by 
Lemma~\ref{lemma5} there exists $\{i_1, i_q\}\in{\cal F}_\Delta$.
Thus by replacing the pair of two sets $S_{j-1}, S_{j}$ by 
an edge $\{i_1, i_q\}$ we 
can remove $S_{j}$ from $Q$. Using the same argument,
we can construct from $Q$ a path $P$ with the required property.

Finally we show the existence of the path.  Assume that if there exist
$i,j\in[n]$ that are not linked by any path, i.e., contained in
different components.  Then considering ${\cal F}_i$ or ${\cal F}_j$ we
immediately know that elements from ${\cal F}_{\Delta}$ contained in
each component must be pairwise disjoint, which entails all the elements
of ${\cal F}_{\Delta}$ are pairwise disjoint. But this contradicts the
condition $(\ref{cond3})$.
\end{proof}

\begin{lemma}
\label{lemma8}
If there exists a length $4$ path
$P = \{\{i_p, i_{p+1}\} \in {\cal F}_{\Delta} \;\vert\; p=1,2,3,4\}$,
then there exists  an edge $\{i_1, i_q\}\in {\cal F}_{\Delta}$ 
with $q=3, 4$ or $5$.
\end{lemma}
\begin{proof}
Consider ${\cal F}_{i_1}$. We have 
$(\{i_3, i_4\}\backslash\{i_1\})\cap (\{i_4, i_5\}\backslash\{i_1\})
=\{i_4\}\ne\emptyset$. Thus by Lemma~\ref{lemma2} and the condition 
$(\ref{cond2})$
there exists $S'\in{\cal F}_{\Delta}$ with $i_1\in S'$ such that 
$(S'\backslash\{i_1\})\subset (\{i_3, i_4\}\backslash\{i_1\}) = \{i_3, i_4\}$
or 
$(S'\backslash\{i_1\})\subset (\{i_4, i_5\}\backslash\{i_1\}) = \{i_4, i_5\}$.
Thus, by taking the condition $(\ref{cond1})$ into account,
we know that $S' = \{i_1, i_q\}$ with $q =3, 4$ or $5$.
\end{proof}

Now we show a converse to what we have proved:
\begin{proposition}
\label{converse}
Consider ${\cal F}_\Delta = \{S_j\;\vert\; S_j\subset [n], j=1,\ldots,\ell\}$ 
satisfying the following conditions, which 
are the same as $(\ref{cond-1})$,  $(\ref{cond0})$ and 
$(\ref{cond1})$:
\begin{enumerate}
\item $S_1\cup\cdots\cup S_\ell = [n]$,
\item $\sharp S_j \geq 2$ for all $1\leq j\leq \ell$, and 
\item $S_i\not\subset S_j$  for all $i,j$ with $i\ne j$.
\end{enumerate}
Then the following are equivalent:
\begin{enumerate}
\item [$(i)$] ${\cal F}_{\Delta}$ satisfies the following conditions,
which are the same as $(\ref{cond2})$ and $(\ref{cond3})$:
\begin{enumerate}
\item 
     $(S\backslash\{i\})\cap (S'\backslash\{i\}) =\emptyset$
     for all distinct $S\backslash\{i\}$ and 
     $S'\backslash\{i\}\in\min{\cal F}_i$ for all $i\in[n]$
\item 
     $S\cap S' \ne\emptyset$
    for some $S, S'\in{\cal F}_{\Delta}$.
\end{enumerate}
\item [$(ii)$] For ${\cal F}_{\Delta}$,
Lemma~\ref{lemma1} 
to Lemma~\ref{lemma8} hold.
\end{enumerate}
\end{proposition}
\begin{proof} 
We have only to show $(i)$ from $(ii)$.
The condition $(\ref{cond3})$ is immediate from 
Lemma~\ref{lemma1}. 
Now we  show the 
condition $(\ref{cond2})$ from Lemmas~\ref{lemma1}$\sim$\ref{lemma8}.

Assume that there exists $i\in [n]$ and $S, T\in{\cal F}_{\Delta}$ such
 that
\begin{equation}
\label{tmp1}
(S\backslash\{i\})\cap (T\backslash\{i\})\ne\emptyset
\end{equation}
and both $S$ and
$T$ are minimal in ${\cal F}_{i}$. 
If $i\in S\cap T$, then $(\ref{tmp1})$ implies $\sharp(S\cap T)\geq 2$,
which contradicts Lemma~\ref{lemma3}. Thus we have either 
$i\notin S$ or $i\notin T$, and we can refine $(\ref{tmp1})$ as 
\begin{equation}
\label{tmp2}
(S\backslash\{i\})\cap (T\backslash\{i\}) = S\cap T = \{j\}
\end{equation}
for some $j$.
Now we have two cases
\begin{description}
\item [case (either $\sharp{S}\geq 3$ or $\sharp{T}\geq 3$)]
We can assume that $\sharp{S}\geq 3$. Then $\sharp{T}=2$ by
Lemma~\ref{lemma4} so that we set $T = \{j,k\}$ for some $k\in [n]$.
If $k=i$, then $i\notin S$ since otherwise we have $T\subset S$,
contradinting the condition $(\ref{cond1})$.
Then we have $\{j\} = (T\backslash\{i\})\subset (S\backslash\{i\}) = S$,
which contradicts 
the assumption that $S$ is  minimal in ${\cal F}_{i}$.
Thus we must have $k\ne i$.

Now assume that $i\in S$. Then by Lemma~\ref{lemma5} we have 
$\{i,k\}\in{\cal F}_\Delta$, so that 
$\{k\}=\{i,k\}\backslash\{i\}\subset T\backslash\{i\} = \{k,j\}$
and $\{i,k\}\ne T$,
which contradicts the minimality of $T$ in ${\cal F}_i$. Thus 
we must have $i\notin S$.
Moreover, if $\in{\cal C}(S)$, then again by Lemma~\ref{lemma5}
we have $\{i,j\}\in{\cal F}_\Delta$ and 
$\{j\}= \{i,j\}\backslash\{i\}\subset T\backslash\{i\}=\{k,j\}$
and $\{i,j\}\ne T$,
which contradicts the minimality of $T$ in ${\cal F}_i$.
Thus we must have $i\notin {\cal C}(S)\cup S$.
Then by
Lemma~\ref{lemma6} we have $\{k,i\}\in {\cal F}_{\Delta}$, so that 
$\{k\}=\{k,i\}\backslash\{i\}\subset T\backslash\{i\}=\{k,j\}$
and $\{k,i\}\ne T$. This contradicts the minimality of 
$T$ in ${\cal F}_{i}$.

\item [case ($\sharp{S}=\sharp{T}=2$)]
We set $S = \{j,h\}$ and $T=\{j,k\}$ for some $h, k (\ne i)\in [n]$.
Then by Lemma~\ref{lemma7} $i$ and $h$ are linked with 
a path 
$P = \{E_1,\ldots, E_r\}$, with edges $E_j\in{\cal F}_{\Delta}$, $j=1,\ldots,
r$, $E_k\cap E_{k+1}\ne\emptyset$ for $k=1,\ldots, r-1$ such that $i\in
E_1$ and $h\in E_r$. If $r\geq 4$, we can take another path
with smaller length by Lemma~\ref{lemma8}, so that we can assume that
$r\leq 3$. Then $P\cup \{S, T\}$ is a path of length $\leq 5$
that links $i$ and $k$. Again by using Lemma~\ref{lemma8} we know that 
$\{i, l\}\in {\cal F}_{\Delta}$ where $l = h, j$ or $k$. But then
we have $\{i,l\}\backslash\{i\}\subset S$
or $\{i,l\}\backslash\{i\}\subset T$ and 
$S$ or $T$ is not minimal in ${\cal F}_{i}$, which contradicts the
assumption.
\end{description}
\end{proof}

\subsection{purity of $\Delta$}
We now characterize the purity of the simplicial complex.
For ${\cal F}_{\Delta}$ as in Prop.~\ref{converse},
we define the associated graph $G_\Delta$ as
follows: The vertex set is $V(G_\Delta) = [n]$ and the edge set is
$E(G_\Delta) = \{ \{i,j\} \;\vert\; \{i,j\}\notin {\cal F}_\Delta\}$.  Then
$G_\Delta$ is exactly 
the skeleton of $\Delta$:
$G_\Delta =  \{F\in\Delta\;\vert\; \dim F\leq
1\}$.

Notice in partcular that, for $S\in{\cal F}_{\Delta}$ with $\sharp{S}\geq 3$,
all the edges in the complete graph ${\cal K}(S)$ over the vertex set
$S$ are contained in $E(G_\Delta)$. 

\begin{remark}
Recall that if $I_\Delta$ is generated in degree $2$, it is called an
{\em edge ideal} represented by a graph $G$ whose edges are suppots of
the generators. In this case,  $G_\Delta$ is nothing but the
complement $\overline{G}$ (see \cite{Vi} Chapter~6). 
\end{remark}

\opn\Simp{Simp}
\opn\Red{Red}

\begin{definition}
For a finite graph $G$ 
we define
\begin{equation*}
  \Simp(G) 
= \left\{
   F\;\vert\;\begin{array}{l}
              \mbox{$\exists H$ a subgraph of $G$ such that}\\
              \mbox{$F=V(H)$ and $H\sim{\cal K}_r$ for some $r$}
	     \end{array}
  \right\}
\end{equation*}
where ${\cal K}_r$ is the complete graph over the 
vertex set $[r]$ and $G\sim H$ denotes isomorphism of 
graphs.
\end{definition}

\begin{definition}
For ${\cal F}_{\Delta}$ as in Prop.~\ref{converse}
and a simplicial complex $\Gamma$ over the vertex set $[n]$,
we define
\begin{equation*}
\Red({\cal F}_\Delta, \Gamma)
= \Gamma\backslash
\{F\in\Gamma
    \;\vert\; 
       S\subset F\mbox{ for some }S\in{\cal F}_\Delta
       \mbox{ with }\sharp{S}\geq 3
\}.
\end{equation*}
\end{definition}

\begin{proposition}[cf. Prop.~6.1.25 \cite{Vi}]
\label{vil6-1.25}
For ${\cal F}_{\Delta}$ as in Prop.~\ref{converse},
we have 
\begin{equation*}
\Delta = \Red({\cal F}_\Delta, \Simp(G_{\Delta})).
\end{equation*}
\end{proposition}
\begin{proof}
$\Delta \subset \Red({\cal F}_\Delta, \Simp(G_{\Delta}))$
is clear from the definitions.
Now take any $F =\{i_1,\ldots, i_r\} 
\in \Red({\cal F}_\Delta, \Simp(G_{\Delta}))$.
Since $F\in\Simp(G_{\Delta})$, 
$\{i_p, i_q\}\in\Delta$, i.e., $\{i_p,i_q\}\notin{\cal F}_\Delta$,
for all $1\leq p<q\leq r$.
In other words, $X_{i_p}X_{i_q}\notin I_\Delta$ for 
all $1\leq p<q\leq r$.
If $F\notin\Delta$, then the monimial $X_{i_1}\cdots X_{i_r}$ is 
in $I_\Delta$.  
Thus there exists a subset $S = \{j_1,\ldots, j_s\}\subset \{i_1,\ldots, i_r\}$
with $3\leq s\leq r$ such that $S\in {\cal F}_\Delta$.
But this contradicts the assumption that 
$F\in\Red({\cal F}_\Delta, \Simp(G_{\Delta}))$.
Thus $F\in\Delta$ as required.
\end{proof}

\begin{corollary}
\label{purity}
For ${\cal F}_{\Delta}$ as in Prop.~\ref{converse},
following are equivalent:
\begin{enumerate}
\item [$(i)$] $\Delta$ is pure.
\item [$(ii)$] for every vertex $i\in G_\Delta$ 
consider a subgraph $H_i\subset G_\Delta$ such that 
\begin{enumerate}
\item $H_i\iso {\cal K}_{r_i}$ and $i\in V(H_i)$
\item $V(H_i)$ contains no subset $S\in{\cal F}_\Delta$
such that $\sharp{S}\geq 3$, and 
\item $r_i$ is maximal among other such subgraphs.
\end{enumerate}
Then $r_i$ is constant for all $i$.
\end{enumerate}
\end{corollary}
\begin{proof}
By Prop.~\ref{vil6-1.25}, we know that a subgraph $H (\subset {G_\Delta})$
isomorphis to ${\cal K}_r$ satisfying the conditions given above 
is exactly the skeleton of an $(r-1)$-face in $\Delta$.
\end{proof}

\subsection{combinatorial characterization of \gci}
Now we obtain a combinatorial characterization of \gci Stanley-Reisner
ring.

\begin{theorem}
\label{theorem:main2}
Let $k[\Delta]$ be a Stanley-Reisner ring with $\Delta = \core\Delta$.
Let $G(I_\Delta) = \{u_1,\ldots, u_\ell\}$ and ${\cal F}_{\Delta} =
\{\supp(u_j)\;\vert\; j=1,\ldots, \ell\}$.  Assume that $k[\Delta]$ is
not a complete intersection.
Then $k[\Delta]$ is a \gci
if and only if the following conditions hold:
\begin{enumerate}
\item for every $S\in{\cal F}_\Delta$  with $\sharp{S}\geq 3$, there exists 
a non-empty set ${\cal C}(S) (\subset [n])$ such that 
   \begin{enumerate}
    \item  ${\cal C}(S)\cap S =\emptyset$,
    \item  for every $i\in {\cal C}(S)$, we have 
	   $E_{ij}:=\{i,j\}\in {\cal F}_\Delta$ for all $j\in S$.
	   Moreover if $S\cap T\ne\emptyset$ for $T\in{\cal F}_\Delta$,
	   then $T = E_{ij}$ for some $i,j$.
    \item  for every $k\notin {\cal C}(S)\cup S$, we have 
	   $\{i,k\}\in {\cal F}_\Delta$ for all $i\in {\cal C}(S)$.
   \end{enumerate}
\item Any two elements $i, j\in [n]$ are linked with a path $P =
\{\{i_k,i_{k+1}\}\;\vert\; k=1,\ldots, r\}$, 
with edges $\{i_k, i_{k+1}\}\in{\cal F}_\Delta$ for  $k=1,\ldots, r$
such that $i=i_1$ and $j = i_{r+1}$.
\item If there exists a length $4$ path $P = \{\{i_p, i_{p+1}\} \in
{\cal F}_\Delta \;\vert\; p=1,2,3,4\}$, then there must be an edge $\{i_1,
i_q\}\in {\cal F}_\Delta$ with $q=3, 4$ or $5$.  
\item for every vertex $i\in G_\Delta$ 
consider a subgraph $H_i\subset G_\Delta$ such that 
\begin{enumerate}
\item $H_i\iso {\cal K}_{r_i}$ and $i\in V(H_i)$
\item $V(H_i)$ contains no subset $S\in{\cal F}_\Delta$
such that $\sharp{S}\geq 3$, and 
\item $r_i$ is maximal among other such subgraphs.
\end{enumerate}
Then $r_i$ is constant for all $i$.
\end{enumerate}
Moreover, if $\Delta\ne\core\Delta$, $k[\Delta]$ is 
a \gci if and only if it is a complete intersection.
\end{theorem}
\begin{proof}
The main part is 
straightforward by Proposition~\ref{converse} and~\ref{vil6-1.25}.
The last part is by Corollary~\ref{ci-is-genci} and 
Lemma~\ref{lemma0}.
\end{proof}

We show a few examples of \gci
Stanely-Reisner ideals.

\begin{example}
Examples of non Cohen-Macaulay edge ideals:
\begin{enumerate}
\item $I_\Delta = (X_1X_3, X_1X_4, X_2X_3, X_2X_4)
=(X_1,X_2)\cap (X_3,X_4)\subset k[X_1,\ldots, X_4]$
and $\Delta = \langle \{1,2\}, \{3,4\}\rangle$ (two disjoint edges).
\item $I_\Delta = (X_1X_2, X_2X_3, X_1X_3,X_3X_4, X_4X_5, X_1X_5)
= (X_2,X_3,X_5)\cap (X_1,X_3,X_5)\cap (X_1,X_3,X_4)\cap (X_1,X_2,X_4)
\subset k[X_1,\ldots, X_5]$ 
and $\Delta = \langle \{1,4\},\{4,2\},\{2,5\},\{5,3\}\rangle$
(a path of length $4$).
\item 
$I_\Delta = (X_1X_2, X_2X_3, X_1X_3, X_3X_4, X_4X_5, X_1X_5, X_2X_5)
= (X_2,X_3,X_5)\cap (X_1,X_3,X_5)\cap  (X_1,X_2,X_4)
\subset k[X_1,\ldots, X_5]$ 
and $\Delta = \langle \{1,4\},\{4,2\}, \{5,3\}\rangle$
(disjoint union a path of length $2$ and an edge).
\item 
$I_\Delta = (X_1X_2, X_1X_5, X_2X_3, X_2X_5, X_3X_4)
= (X_2,X_4,X_5)\cap (X_2,X_3,X_5)\cap (X_1,X_2,X_4)\cap (X_1,X_2,X_3)
\cap (X_1,X_3,X_5)\subset k[X_1,\ldots, X_5]$
and\\
$\Delta = \langle \{1,3\}, \{3,5\},\{5,4\}, \{4,1\} \{4,2\}\rangle$
(an edge attached to a circle).
\item 
$I_\Delta = (X_1,\ldots, X_n)\cap (X_{n+1},\ldots, X_{2n})
\subset k[X_1,\ldots, X_{2n}]$. Notice that $G(I_\Delta)$ is 
a bipartite graph.
\end{enumerate}
\end{example}

\begin{example}
Examples of Cohen-Macaulay edge ideals:
\begin{enumerate}
\item $I_\Delta = 
(X_1X_2, X_2X_3, X_3X_4) = (X_2,X_4)\cap (X_2,X_3)\cap (X_1,X_3)
\subset k[X_1,\ldots, X_4]$ 
and $\Delta = \langle \{1,3\}, \{1,4\}, \{4,2\}\rangle$ (a path of
      length $3$)
\item $I_\Delta = (X_1X_2, X_2X_3, X_3X_4, X_4X_5, X_5X_1)
= (X_2,X_4,X_5)\cap (X_1,X_2,X_4)\cap (X_1,X_3,X_4)\cap (X_1,X_3,X_5)
\cap (X_2,X_3,X_5)\subset k[X_1,\ldots, X_5]$ 
and\\
$\Delta = \langle \{1,3\}, \{3,5\}, \{5,2\}, \{2,4\}, \{4,1\}\rangle$
(a circle)
\end{enumerate}
For Cohen-Macaulay edge ideals, see \cite{Vi}.
\end{example}

\begin{example}
Ideals whose generators contain degree $\geq 3$ monomials:
\begin{enumerate}
\item $I_\Delta = (X_1X_2X_3) + (X_1, X_2, X_3)*(X_4,X_5, X_6) +
(X_4X_7, X_5X_7, X_6X_7)\subset k[X_1,\ldots, X_7]$
and 
\begin{equation*}
\Delta = \langle
           \{1,2,7\}, \{1,3,7\}, \{2,3,7\},
               \{4,5,6\}
         \rangle,
\end{equation*}
which is a not Cohen-Macaulay complex since it is not connected.

\item $I_\Delta = (X_1X_2X_3X_4) 
+ (X_1,X_2,X_3,X_4)*(X_5,X_6,X_7)$ 
and 
\begin{equation*}
\Delta = \langle
           \{1,2,3\},\{1,2,4\},\{1,3,4\},\{2,3,4\},
           \{5,6,7\}
         \rangle,
\end{equation*}
which is a not Cohen-Macaulay complex since it is not connected.
\end{enumerate}
\end{example}

\bibliographystyle{amsplain} 

\begin{thebibliography}{99}

\bibitem{AV}
R.\ Achilles and W.\ Vogel, \"Uber vollst\"andige 
Durchschnitte in lokalen Ringen. Math. Nachr. 89 (1979), 285--298. 

\bibitem{BH} W.\ Bruns and J.\ Herzog, {\em Cohen-Macaulay Rings},
Revised version, Cambridge Studies in Advanced Mathematics,~39.
Cambridge University Press, Cambridge, 1998.

\bibitem{CN} R.\ C.\ Cowsik and M.\ V.\ Nori,
On the fibres of blowing up. J. Indian Math. Soc. (N.S.) 40 
(1976), no. 1-4, 217--222 (1977). 

\bibitem{HO}M.\ Herrmann and U.\ Orbanz,
Faserdimensionen von Aufblasungen lokaler Ringe und 
\"Aquimultiplizitat. J. Math. Kyoto Univ. 20 (1980), no. 4, 651--659. 

\bibitem{S} R.\.P.\ Stanley, {\em Combinatorics and
commutative algebra}, 
Second edition. Progress in Mathematics,~41. 
Birkh\"auser Boston, Inc., Boston, MA, 1996.

\bibitem{T} Y.\ Takayama, Combinatorial characterizations 
of generalized Cohen-Macaulay monomial ideals,
Bull. Math. Soc. Sci. Math. Roumanie (N.S.) 48(96) (2005), no. 3, 327--344

\bibitem{Vi} R.\ H.\ Villarreal, {\em Monomial Algebras},
Monographs and Textbooks in Pure and Applied Mathematics,~238. 
Marcel Dekker, Inc., New York, 2001.

\bibitem{Wal} R.\ Waldi, Vollst\"{a}ndige Durchschnitte in Cohen-Macaulay-Ringen. 
Arch. Math. (Basel) 31 (1978/79), no. 5, 439--442

\end{thebibliography}


\end{document}